\documentclass[12pt,reqno]{article}

\usepackage[usenames]{color}
\usepackage{amssymb}
\usepackage{graphicx}
\usepackage{amscd}

\usepackage[colorlinks=true,
linkcolor=webgreen,
filecolor=webbrown,
citecolor=webgreen]{hyperref}

\definecolor{webgreen}{rgb}{0,.5,0}
\definecolor{webbrown}{rgb}{.6,0,0}

\usepackage{color}
\usepackage{fullpage}
\usepackage{float}

\usepackage{graphics,amsmath,amssymb}
\usepackage{amsthm}
\usepackage{amsfonts}
\usepackage{latexsym}
\usepackage{epsf}

\newcommand{\seqnum}[1]{\href{http://oeis.org/classic/#1}{\underline{#1}}}

\begin{document}

\newtheorem{theorem}{Theorem}[section]
\newtheorem{proposition}{Proposition}[section]
\newtheorem{corollary}{Corollary}[section]
\newtheorem{lemma}{Lemma}[section] 

\newtheorem{notation}{Notation}
\newtheorem{remark}{Remark}
\newtheorem{example}{Example}
\newtheorem{conjecture}{Conjecture}
\newcommand{\p}{\textit{Proof. }}
\newcommand{\q}{\hfill $\Box$}
\newcommand{\lt}{\left( }
\newcommand{\rt}{\right) }

\begin{center}
\vskip 1cm{\LARGE\bf 
Decomposition into {\em weight $\times$ level $+$ jump}\\
and application to a new classification of primes}
\vskip 1cm
\large
R\'emi Eismann\\
\href{mailto:reismann@free.fr}{\tt reismann@free.fr}\\
\end{center}

\vskip .2 in

\begin{abstract}
In this paper we introduce an Euclidean decomposition of elements $a_n$ of an increasing sequence of natural numbers $(a_n)_{n\in\mathbb{N}^*}$ 
into {\em weight $\times$ level $+$ jump} which we use to classify the numbers $a_n$ either by {\em weight} or by {\em level}. We then show 
that this decomposition can be seen as a generalization of the sieve of Eratosthenes (which is the particular case of the whole sequence of natural numbers). 
We apply this decomposition to prime numbers in order to obtain a new classification of primes, we analyze a few properties of 
this classification and we make a series of conjectures based on numerical data. Finally we show how composite numbers and $2-$almost primes behave under the decomposition.

\end{abstract}

\section{Decomposition algorithm of numbers into {\em weight $\times$ level $+$ jump} and application to a classification scheme}
\label{s:1}
We introduce an algorithm whose input is an increasing sequence of positive integers $(a_n)_{n\in\mathbb{N}^*}$ 
and whose output is a sequence of unique triplets of positive integers $(k_n,L_n,d_n)_{n\in\mathbb{N}^*}$.
\bigskip

We define the {\em jump} (first difference, gap) of $a_n$ by
\[ 
d_n:=a_{n+1}-a_n\text{.}
\] 

Then let $l_n$ be defined by 
\[ 
l_n:=\left\{
\begin{array}{l} 
a_n-d_n \text{ if } a_n-d_n>d_n \text{;}\\ 
0 \text{ otherwise. }
\end{array}
\right.
\] 

The {\em weight} of $a_n$ is defined to be 
\[ 
k_n:=\left\{
\begin{array}{l} 
\min\{k\in\mathbb{N}^* \text{ s.t. } k>d_n, k|l_n \} \text{ if } l_n\neq 0 \text{;}\\ 
0 \text{ otherwise. }
\end{array}
\right.
\] 

Finally we define the {\em level} of $a_n$ by 
\[
L_n:=\left\{
\begin{array}{l} 
\frac{l_n}{k_n} \text{ if } k_n\neq 0 \text{;}\\
0 \text{ otherwise. }
\end{array}
\right.
\]

We then have a decomposition of $a_n$ into {\em weight $\times$ level $+$ jump}: $a_n=l_n+d_n=k_n \times L_n + d_n$ when $l_n\neq 0$. 

In the Euclidean division of $a_n$ by its {\em weight} $k_n$, the quotient is the {\em level} $L_n$, and the remainder is the {\em jump} $d_n$. 

\begin{lemma} 
\label{t:1}
A necessary and sufficient condition for the decomposition of a number $a_n$ belonging to an increasing sequence of positive 
integers $(a_n)_{n\in\mathbb{N}^*}$ into weight $\times$ level $+$ jump to hold is that 
\[
\begin{array}{cl}
 a_{n+1}<\frac{3}{2}a_n\text{.} \\
\end{array}
\].
\end{lemma}

\begin{proof} 
The decomposition is possible if $l_n\neq0$, that is if $a_n-d_n > d_n$, which can be rewritten as $a_{n+1}<\frac{3}{2}a_n$.
\end{proof}

In order to use this algorithm to classify the numbers $a_n$ we introduce the following rule (whose meaning will become clearer 
in the next section):   if for $a_n$ we have $k_n=L_n=l_n=0$ then $a_n$ is not classified; if for $a_n$ we have $k_n>L_n$ then $a_n$ is classified by {\em level}, if not then $a_n$ is classified by {\em weight}.

\section{Application of the algorithm to the sequence of natural numbers}
\label{s:2}

In this situation we have $a_n=n$ et $d_n=1$. The decomposition is impossible for $n=1$ and $n=2$ ($l_1 = l_2 = 0$).
Apart from those two cases, we have the decomposition of $n$ into {\em weight $\times$ level $+$ jump}: 
$n=k_n \times L_n+1$ when $n > 2$ and we also have the following relations 
\[
\begin{array}{cl}
 & L_n=1 \\
\Leftrightarrow & k_n>L_n\\
\Leftrightarrow & k_n=l_n=n-1\\
\Leftrightarrow & l_n=n-1 \text{ is prime, }\\
\end{array}
\]

\[
\begin{array}{cl}
 & L_n>1 \\
\Leftrightarrow & k_n \leq L_n\\
\Leftrightarrow & k_n \times L_n=l_n=n-1 \\
\Leftrightarrow & l_n=n-1 \text{ is composite. }\\
\end{array}
\]

The {\em weight} of $n$ is the smallest prime factor of $n - 1$ and the {\em level} of $n$  is the largest proper divisor of $n - 1$ .
We can characterize the fact that a number $l_n=n-1$ is prime by the fact that $n$ is classified by {\em level} 
(or equivalently here by the fact that $n$ is of {\em level} $1$).
 
Since there is an infinity of prime numbers, there is an infinity of natural numbers of {\em level} $1$. Similarly there is an infinity of natural numbers with a {\em weight } equal to $k$ with $k$ prime. The algorithm allows to separate prime numbers ($l_n$ or {\em weights } of natural numbers classified by {\em level}) from composite numbers ($l_n$ of natural numbers classified by {\em weight}), and is then indeed 
a reformulation of the sieve of Eratosthenes. Apply this algorithm to any other increasing sequence 
of positive integers, for example to the sequence of prime numbers itself, can be seen as an generalization of that sieve. 

\begin{table}[th!]
\begin{center}
\begin{tabular}{|c|c|c|c|c|} 
\hline $n$ \seqnum{A000027} & $k_n$ \seqnum{A020639}($n-1$) & $L_n$ \seqnum{A032742}($n-1$) & $d_n$ & $l_n$\\
\hline 1&0&0&1&0 \\
\hline 2&0&0&1&0 \\
\hline 3&2&1&1&2 \\
\hline 4&3&1&1&3 \\
\hline 5&2&2&1&4 \\
\hline 6&5&1&1&5 \\
\hline 7&2&3&1&6 \\
\hline 8&7&1&1&7 \\
\hline 9&2&4&1&8 \\
\hline 10&3&3&1&9 \\
\hline 11&2&5&1&10 \\
\hline 12&11&1&1&11 \\
\hline 13&2&6&1&12 \\
\hline
\end{tabular}
\end{center}
\caption{ The 13 first terms of the sequences of {\em weights}, {\em levels}, {\em jumps}, and $l_n$ in the case of the sequence of natural numbers.}
\label{tab2} 
\end{table}

\begin{figure}[th!]
\begin{center}
\input{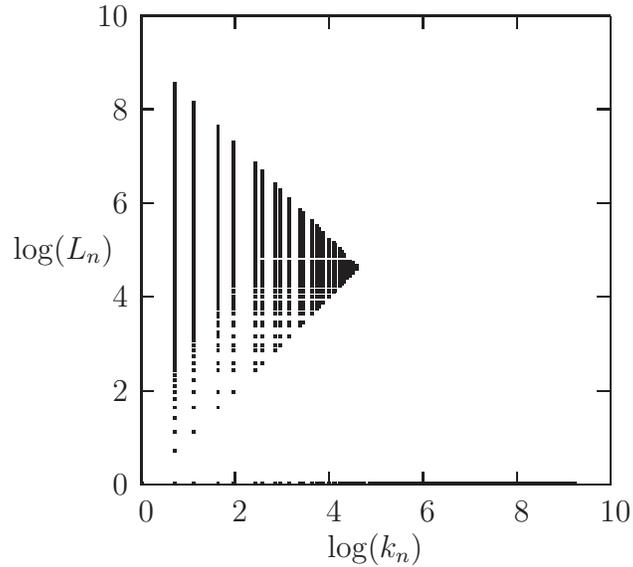}
\caption{Plot of natural numbers in $\log(k_n)$ vs. $\log(L_n)$ coordinates (with $n\leq10000$). The sieve of Eratosthenes.} 
\label{plot_entiers}
\end{center}
\end{figure}

\newpage
\section{Application of the algorithm to the sequence of primes}
\label{s:3}

We can wonder what happens if we try to apply the decomposition to the sequence of primes itself: 
for any $n\in\mathbb{N}^*$ we have $a_n=p_n$ and $d_n=g_n$ (the prime gap). 
The algorithm of section~\ref{s:1} can then be rewritten with these new notations as follows. 

The {\em jump} (first difference, gap) of $p_n$ is 
\[ 
g_n:=p_{n+1}-p_n\text{.}
\] 

Let $l_n$ be defined by
\[ 
l_n:=\left\{
\begin{array}{l} 
p_n-g_n \text{ if } p_n-g_n>g_n \text{;}\\ 
0 \text{ otherwise.}
\end{array}
\right.
\] 

The {\em weight} of $p_n$ is then  
\[ 
k_n:=\left\{
\begin{array}{l} 
\min\{k\in\mathbb{N}^* \text{ s.t. } k>g_n, k|l_n \} \text{ if } l_n\neq 0 \text{;}\\ 
0 \text{ otherwise. }
\end{array}
\right.
\] 

The {\em level} of $p_n$ is 
\[
L_n:=\left\{
\begin{array}{l} 
\frac{l_n}{k_n} \text{ if } k_n\neq 0 \text{;}\\
0 \text{ otherwise. }
\end{array}
\right.
\]

In the Euclidean division of $p_n$ by its {\em weight} $k_n$, the quotient is the {\em level} $L_n$, and the remainder is the {\em jump} $g_n$. 

So we have the decomposition of $p_n$ into {\em weight $\times$ level $+$ jump} reads $p_n=k_n \times L_n+g_n$ when $l_n \neq 0$. So one 
should investigate for which $n$ we have $l_n \neq 0$, which is provided by the following result.

\begin{theorem} 
\label{t:2}
This decomposition is always possible except for $p_1=2$, $p_2=3$ and $p_4=7$ \\
(i.e., $p_{n+1}\geq \frac{3}{2}p_n$ holds only for $n=1$, $n=2$ and $n=4$).
\end{theorem}

\begin{proof} 
The decomposition if possible if, and only if, $l_n$ is not equal to zero. 
But $l_n\neq0$ if and only if $p_{n+1}<\frac{3}{2}p_n$, that is $p_n-g_n>g_n (*)$ by lemma~\ref{t:1}.
Let us now apply results of Pierre Dusart on the prime counting function $\pi$ to show that this is always true except 
for $n=1$, $n=2$ and $n=4$. 

Indeed this last equation $(*)$ can be rewritten in terms of $\pi$ as $\pi(\frac{3}{2}x)-\pi(x)>1$ 
(i.e., there is always a number strictly included between $x$ and $\frac{3}{2}x$ for any $x\in\mathbb{R}^+$). 
But Dusart has shown \cite{pD01,pD02} that on the one hand for $x\geq599$ we have  
\[
\pi(x)\geq \frac{x}{\log x}\Big(1+\frac{1}{\log x}\Big)
\]

and on the other hand for $x>1$ we have 

\[
\pi(x)\leq \frac{x}{\log x}\Big(1+\frac{1.2762}{\log x}\Big)
\]

So for $x\geq 600$ we have 

\[
\pi(\frac{3}{2}x)-\pi(x) > \frac{900}{\log 900}\Big(1+\frac{1}{\log 900}\Big) 
- \frac{600}{\log 600}\Big(1+\frac{1.2762}{\log 600}\Big) 
\]
and since the right hand side of this inequality is approximately equal to $39.2$ we indeed have that 
$\pi(\frac{3}{2}x)-\pi(x)>1$, so the inequality $(*)$ holds for any prime greater than $600$. 
We check numerically that it also holds in the remaining cases when $x<600$, except for the aforementioned exceptions 
$n=1$, $n=2$ and $n=4$ which ends the proof.
\end{proof}

Let us now state a few direct results. For any $p_n$ different from $p_1=2$, $p_2=3$ and $p_4=7$ we have 
\[
\begin{array}{cl}
 \gcd(g_n,2)=2\text{, }\\
 \gcd(p_n,g_n) = \gcd(p_n - g_n,g_n) = \gcd(l_n,g_n) = \gcd(L_n,g_n) = \gcd(k_n,g_n) = 1\text{, }\\
 3 \le k_n \le l_n\text{, }\\
 1 \le L_n \le \frac{l_n}{3}\text{, }\\
 2 \le g_n \le k_n - 1\text{, }\\
 2 \times g_n + 1 \le p_n\text{.}\\
\end{array}
\]

\begin{lemma} 
\label{t:3}
$p$ is a prime such that $p>3$ and $p+2$ is also prime if and only if $p$ has a weight equal to $3$.
\end{lemma}

\begin{proof} 
The primes $p>3$ such that $p+2$ is also prime are of the form $6n-1$, so $p-2$ is of the form $6n-3$.  
The smallest divisor greater than $2$ of a number of the form $6n-3$ is $3$. If $p_n$ has a weight equal 
to $3$ then $p_n>3$ and the jump $g_n$ is equal to $2$ since we know that $2 \le g_n \le k_n - 1$ and 
$2 \times g_n + 1 \le p_n$.
\end{proof}

\begin{table}[th!]
\begin{center}
\begin{tabular}{|c|c|c|c|c|c|} 
\hline $n$ & $p_n$ \seqnum{A000040} & $k_n$ \seqnum{A117078} & $L_n$ \seqnum{A117563} & $d_n$ \seqnum{A001223} & $l_n$ \seqnum{A118534}\\
\hline 1&2&0&0&1&0 \\
\hline 2&3&0&0&2&0 \\
\hline 3&5&3&1&2&3 \\
\hline 4&7&0&0&4&0 \\
\hline 5&11&3&3&2&9 \\
\hline 6&13&9&1&4&9 \\
\hline 7&17&3&5&2&15 \\
\hline 8&19&5&3&4&15 \\
\hline 9&23&17&1&6&17 \\
\hline 10&29&3&9&2&27 \\
\hline 11&31&25&1&6&25 \\
\hline 12&37&11&3&4&33 \\
\hline 13&41&3&13&2&39 \\
\hline 14&43&13&3&4&39 \\
\hline 15&47&41&1&6&41 \\
\hline 16&53&47&1&6&47 \\
\hline 17&59&3&19&2&57 \\
\hline
\end{tabular}
\end{center}
\caption{ The 17 first terms of the sequences of {\em weights}, {\em levels}, {\em jumps} and $l_n$ in the case of the sequence of primes.}
\label{tab3} 
\end{table}

\begin {figure}[th!]
\begin{center}
\input{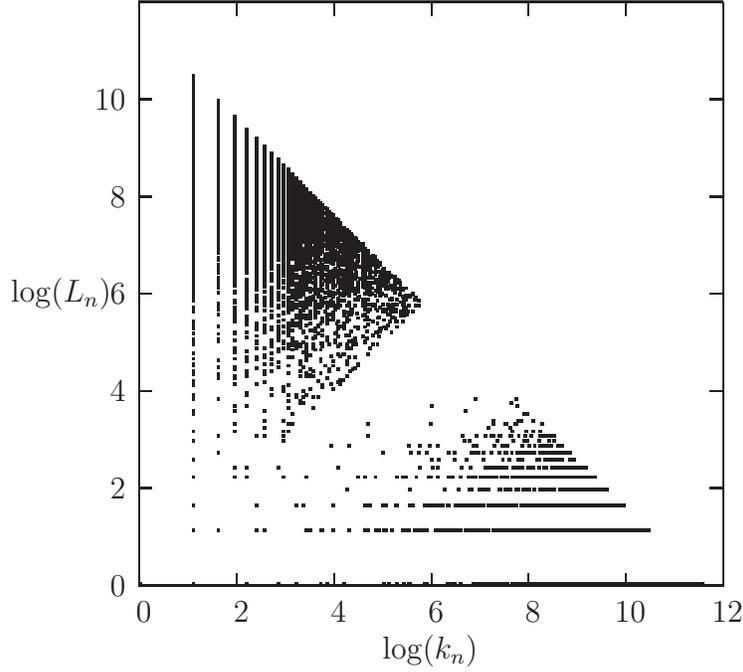}
\caption{Plot of prime numbers in $\log(k_n)$ vs. $\log(L_n)$ coordinates (with $n\leq10000$).} 
\label{plot_premiers}
\end{center}
\end {figure}

\newpage
\section{Classification of prime numbers}
\label{s:4}

We introduce the following classification principle: \\
- if for $p_n$ we have $k_n=L_n=l_n=0$ then $p_n$ is not classified;\\
- if for $p_n$ we have $k_n>L_n$ then $p_n$ is classified by {\em level}, if not $p_n$ is classified by {\em weight}; \\
- furthermore if for $p_n$ we have that $l_n$ is equal to some prime $p_{n-i}$ then $p_n$ is of {\em level} $(1;i)$.\\

For $n\leq5.10^7$, $17,11 \%$ of the primes $p_n$ are classified by {\em level} and $82,89 \%$ are classified by {\em weight}. 
\bigskip

We have the following direct results:

If $p_n$ is classified by {\em weight} then
\[
\begin{array}{cl}
 g_n + 1 \le k_n \le \sqrt{l_n} \le L_n \le \frac{l_n}{3}\text{.}\\
\end{array}
\]

\begin{table}[th!]
\begin{center}
\begin{tabular}{|c|c|c|c|} 
\hline Number of & which have & \% / total of primes & \% / total\\
primes & a {\em weight} equal to  & classified by {\em weight} &  \\
\hline 3370444&3&8.132&6.741\\
\hline 1123714&5&2.711&2.247\\
\hline 1609767&7&3.884&3.219\\
\hline 1483560&9&3.579&2.967\\
\hline 1219514&11&2.942&2.439\\
\hline 1275245&13&3.077&2.550\\
\hline 1260814&15&3.042&2.522\\
\hline 1048725&17&2.530&2.097\\
\hline 1051440&19&2.546&2.103\\
\hline 1402876&21&3.385&2.806\\
\hline 893244&23&2.155&1.786\\
\hline
\end{tabular}
\end{center}
\caption{Distribution of primes for the 11 smallest {\em weights} (with $n\leq5.10^7$).}
\label{tab4} 
\end{table}

If $p_n$ is classified by {\em level} then 
\[
\begin{array}{cl}
 L_n < \sqrt{l_n} < k_n \le l_n\text{;}\\
L_n + 2 \le g_n + 1 \le k_n \le l_n\text{.}\\
\end{array}
\]

Pimes $p_n$ for which $g_n > \sqrt(l_n)$ are : $2$, $3$, $5$, $7$, $13$, $19$, $23$, $31$, $113$ for  $n\leq5.10^7$.

\begin{table}[th!]
\begin{center}
\begin{tabular}{|c|c|c|c|} 
\hline Number of & which are & \% / total of primes & \% / total\\
primes & of {\em level}  & classified by {\em level} &  \\
\hline 2664810&1&31.15&5.330\\
\hline 2271894&3&26.56&4.544\\
\hline 963665&5&11.27&1.927\\
\hline 444506&7&5.197&0.8890\\
\hline 640929&9&7.493&1.282\\
\hline 254686&11&2.978&0.5094\\
\hline 155583&13&1.819&0.3112\\
\hline 351588&15&4.110&0.7032\\
\hline 115961&17&1.356&0.2319\\
\hline 78163&19&0.9138&0.1563\\
\hline 148285&21&1.734&0.297\\
\hline
\end{tabular}
\end{center}
\caption{Distribution of primes for the 11 smallest {\em levels} (with $n\leq5.10^7$).}
\label{tab5} 
\end{table}

\newpage
If $p_n$ is of {\em level} $(1;i)$ then  
\[
\begin{array}{cl}
    L_n = 1 \text{ and } l_n = k_n = p_{n-i}\text{,}\\
    p_n = p_{n-i} + g_n \text{ or } p_{n+1} - p_n = p_n - p_{n-i}\text{.}\\
\end{array}
\]

Primes of {\em level} $(1;1)$ are the so-called "balanced primes" (\seqnum{A006562}). If $p_n$ is of {\em level} $(1;1)$ then  
\[
\begin{array}{cl}
    L_n = 1 \text{ and } l_n = k_n = p_{n-1}\text{,}\\
    g_n = g_{n-1} \text{ or } p_{n+1} - p_n = p_n - p_{n-1}\text{,}\\
    p_n = \frac{p_{n+1} + p_{n-1}} {2}\text{.}\\
\end{array}
\]

\begin{table}[th!]
\begin{center}
\begin{tabular}{|c|c|} 
\hline Number of & which are of {\em level} $(1;i)$\\
primes & $i$ \\
\hline 1307356&1\\
\hline 746381&2\\
\hline 345506&3\\
\hline 153537&4\\
\hline 65497&5\\
\hline 27288&6\\
\hline 11313&7\\
\hline
\end{tabular}
\end{center}
\caption{Distribution of primes of {\em level} $(1;i)$ (with $n\leq5.10^7$, $i\leq7$).}
\label{tab6} 
\end{table}

\begin{table}[th!]
\begin{center}
\begin{tabular}{|c|c|c|c|c|} 
\hline n $<=$ & number of primes & number of primes & \% primes & \% primes \\
& classified by level & classified by weight & classified by level & classified by weight \\
\hline 100&44&53&44&53\\
\hline 1000&324&673&32,4&67,3\\
\hline 10000&2766&7231&27,66&72,31\\
\hline 100000&22999&76998&23&77\\
\hline 1000000&203441&796556&20,34&79,66\\
\hline 10000000&1828757&8171240&18,29&81,71\\
\hline 50000000&8553468&41446529&17,11&82,89\\
\hline
\end{tabular}
\end{center}
\caption{repartition between primes classified by {\em level} and primes classified by {\em weight}.}
\label{tab7} 
\end{table}

\newpage
\section{Conjectures on primes}
\label{s:5}

From our numerical data on the decomposition of primes $p_n$ until $n=5.10^7$ we make the following conjectures.

Since we have shown previously that the smallest number of each twin prime pair (except $3$) has a {\em weight} equal to $3$,
 the well-known conjecture on the existence of an infinity of twin primes can be rewritten as
\begin{conjecture} 
\label{c:1}
The number of primes with a weight equal to $3$ is infinite.
\end{conjecture}

To extend this conjecture, and by analogy with the decomposition of natural numbers for which we know that for any prime $k$ there exist 
an infinity of natural numbers with a {\em weight} equal to $k$ and that there exist an infinity of natural numbers of {\em level} $1$, 
we make this two conjectures
\begin{conjecture} 
\label{c:2}
The number of primes with a weight equal to $k$ is infinite for any $k\ge 3$ which is not a multiple of $2$.
\end{conjecture}

\begin{conjecture} 
\label{c:3}
The number of primes of level $L$ is infinite for any $L\ge 1$ which is not a multiple of $2$.
\end{conjecture}

Now, based on our numerical data and again by analogy with the decomposition of natural numbers for which we know that the natural 
numbers which are classified by {\em level} have a $l_n$ or a {\em weight} which is always prime we conjecture
\begin{conjecture} 
\label{c:4}
Except for $p_6=13$, $p_{11}=31$, $p_{30}=113$, $p_{32}=131$ et $p_{154}=887$, primes which are classified by level have a 
weight which is itself a prime.
\end{conjecture}

The conjecture on the existence of an infinity of balanced primes can be rewritten as 
\begin{conjecture} 
\label{c:5}
The number of primes of level $(1;1)$ is infinite.
\end{conjecture}
That we can easily generalize by
\begin{conjecture} 
\label{c:6}
The number of primes of level $(1;i)$ is infinite for any $i\ge 1$.
\end{conjecture}

We make the following conjectures, for which we have no rigorous arguments yet
\begin{conjecture} 
\label{c:7}
If the jump $g_n$ is not a multiple of $6$ then $l_n$ is a multiple of $3$.
\end{conjecture}

\begin{conjecture} 
\label{c:8}
If $l_n$ is not a multiple of $3$ then jump the $g_n$ is a multiple of $6$.
\end{conjecture}

According to the numerical data and knowing that the primes are rarefying among the natural numbers, we make the following conjecture: 
\begin{conjecture} 
\label{c:9}
The primes classified by level are rarefying among the prime numbers.
\end{conjecture}

\newpage 
\section{Decomposition of composite numbers and of $2-$almost primes.}
\label{s:6}

In this section we only provide the plots of the distribution of composite numbers and of $2-$almost primes in $\log(k_n)$ vs. $\log(L_n)$ coordinates.

\begin {figure}[th!]
\begin{center}
\input{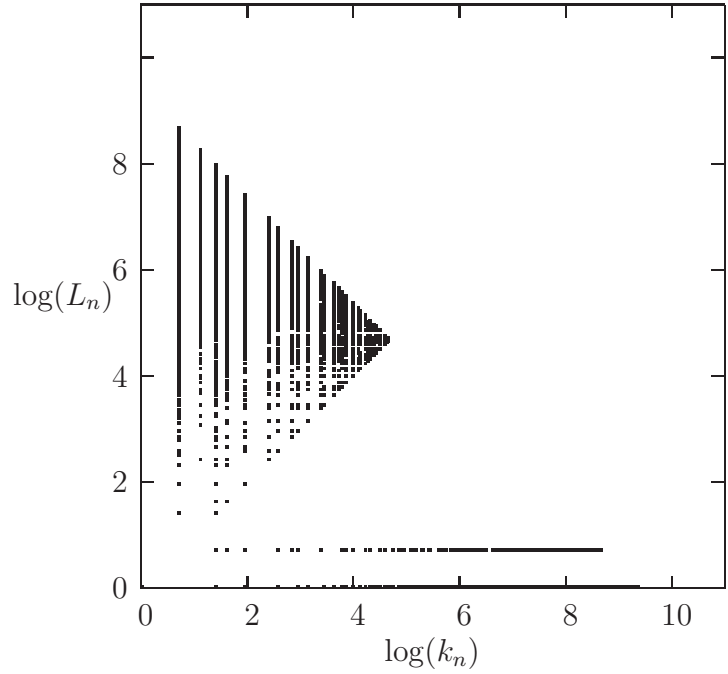}
\caption{Plot of composite numbers (\seqnum{A002808}) in $\log(k_n)$ vs. $\log(L_n)$ coordinates (with $n\leq9999$).}
\label{plot_composite}
\end{center}
\end {figure}
The sequence of {\em weights} of composite numbers is \seqnum{A130882}.
 
\newpage
\begin {figure}[th!]
\begin{center}
\input{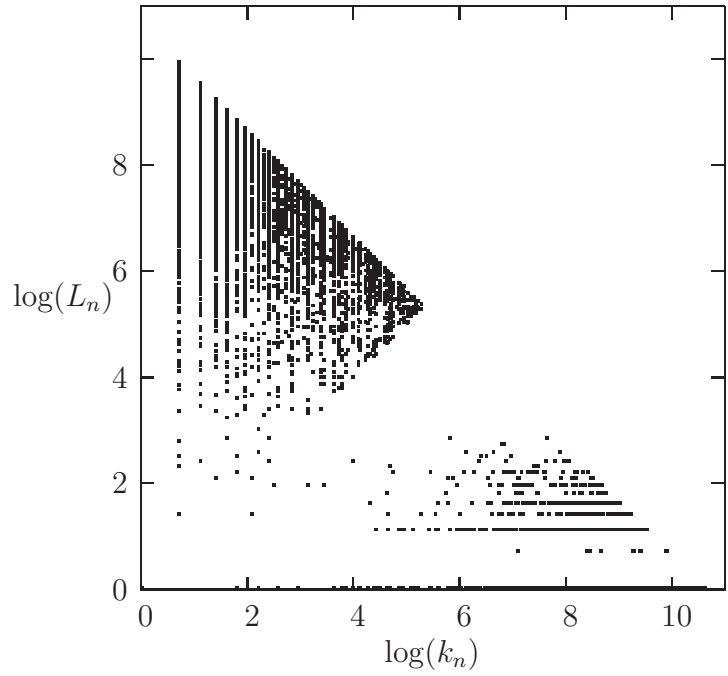}
\caption{Plot of $2-$almost primes (\seqnum{A001358}) in $\log(k_n)$ vs. $\log(L_n)$ coordinates (with $n\leq9999$).}
\label{plot_semi}
\end{center}
\end {figure}
The sequence of {\em weights} of $2-$almost primes is \seqnum{A130533}.

\newpage
\section{Acknowledgements}
The author wishes to thank Thomas Sauvaget for the proof of theorem~\ref{t:2}, his advices and translating the text into English.
The author whishes to thank Fabien Sibenaler for developping programs in Java and Assembly implementing the decomposition algorithm 
and for his encouragements. The author also whishes to thank Jean-Paul Allouche for his advices and for his encouragements and 
N.~J.~A.~Sloane for its help through the OEIS.

\bigskip
\hrule
\bigskip

\noindent 2000 {\it Mathematics Subject Classification}: 
Primary 11B83; Secondary 11P99.

\noindent \emph{Keywords: } 
integer sequences, Eratosthenes, sieve, primes, classification, prime gaps, twin primes, balanced primes.

\bigskip
\hrule
\bigskip

\noindent (Concerned with sequences
\seqnum{A000027}, 
 \seqnum{A020639},  
 \seqnum{A032742},  
 \seqnum{A000040},  
 \seqnum{A117078},  
 \seqnum{A117563},  
 \seqnum{A001223},  
 \seqnum{A118534},  
 \seqnum{A001359},  
 \seqnum{A006562},  
 \seqnum{A125830},  
 \seqnum{A117876}, 
 \seqnum{A074822}, 
 \seqnum{A002808},
 \seqnum{A130882},
 \seqnum{A001358} and 
 \seqnum{A130533}.)

\bigskip
\hrule
\bigskip

\end{document}